\documentclass[reqno,11pt]{amsart}
\usepackage{epsfig}\usepackage{amsfonts,amsmath,amssymb,mathrsfs,verbatim,amsthm,amscd}
\numberwithin{equation}{section}
\newcommand{\beq}{\begin{equation}}
\newcommand{\ee}{\end{equation}}
\newcommand{\bea}{\begin{eqnarray}}
\newcommand{\eea}{\end{eqnarray}}

\usepackage{hyperref}
\def\stackreb#1#2{\ \mathrel{\mathop{#1}\limits_{#2}}}

\newcommand{\CC}{\mathbb C}
\newcommand{\RR}{\mathbb R}
\newcommand{\Z}{\mathbb Z}

\newcommand{\T}{\mathbb T}

\begin{document}

\vspace*{-2em}

\title[Elliptic hypergeometric functions: integrals versus series]
{Elliptic hypergeometric functions: \\ integrals versus series}

\medskip

\author{Vyacheslav\, P.~ Spiridonov}%

\address{Laboratory of Theoretical Physics, JINR, Dubna, Moscow region, Russia and
National Research University Higher School of Economics, Moscow, Russia
}

\vspace*{-1.0em}

\begin{abstract}
The univariate elliptic beta integral is represented as a bilinear combination of
infinite $_{10}V_9$ very-well-poised elliptic hypergeometric series representing
the sum of residues of the integrand poles. Convergence of this combination of series
for some particular choice of parameters is discussed.
Additionally, the asymptotics of the Frenkel--Turaev sum for a terminating
$_{10}V_9$ series is considered when the termination parameter $n$ goes to infinity.
\end{abstract}


\maketitle

\tableofcontents

\section{Introduction}

Special functions \cite{aar} form a thin layer of universal mathematical objects with ubiquitous variety
of applications in natural sciences.
Elliptic hypergeometric functions \cite{spi:aar} currently lie on the top of special functions
of hypergeometric type. They appeared first in the form of terminating elliptic hypergeometric
series depending on the basic variable $q$, elliptic nome $p$ and some
number of free abelian parameters  \cite{FT}. For a special set of abelian parameters
they emerged in the hidden form as elliptic function solutions of the
Yang--Baxter equation of RSOS type \cite{djkmo}. The elliptic hypergeometric
functions which are truly transcendental over the field of elliptic functions were
discovered in \cite{spi:umn} in the form of elliptic hypergeometric integrals
depending on two basic variables $q$ and $p$ on equal footing.

Until recently it was not known whether the infinite (non-terminating) elliptic hypergeometric
series converge for some choice of parameters or not. This question was partially clarified
in \cite{KS}, where a nontrivial choice of parameters was suggested for
which the infinite series converge and define double periodic functions of parameters with natural
boundaries inside fundamental parallelograms of periods. In the present paper
we continue this research and represent the univariate elliptic beta integral \cite{spi:umn}
as a bilinear combination of infinite $_{10}V_9$ very-well-poised elliptic hypergeometric
series through the residue calculus. Convergence of this combination of series
for some particular choice of the basic variables $p$ and $q$ is discussed. However,
it is not established yet.
The criterion of convergence of infinite elliptic hypergeometric series described
in \cite{KS} is also applied to the investigation of asymptotics of the Frenkel--Turaev
sum for a terminating $_{10}V_9$ series \cite{FT} when the termination parameter $n$
goes to infinity.

\section{Elliptic beta integral as a sum of residues}

The elliptic beta integral is a univariate contour integral of hypergeometric type
admitting exact evaluation \cite{spi:umn}
\begin{eqnarray}\nonumber && \makebox[-0em]{}
I(\underline{t};p,q):=\frac{(p;p)_\infty(q;q)_\infty}{4\pi  i }
\int_\T\frac{\prod_{a=1}^6
\Gamma(t_az;p,q)\Gamma(t_az^{{- 1}};p,q)}{\Gamma(z^2;p,q)\Gamma(z^{-2};p,q)}\frac{dz}{z}
\\ && \makebox[2em]{}
=\prod_{1\leq a<b\leq6}\Gamma(t_at_b;p,q), \qquad  \prod_{j=1}^6t_a=pq,
\label{ellbeta}\end{eqnarray}
where $\T$ is the unit circle of positive orientation, $|t_a|<1$, $(z;p)_\infty=\prod_{j=1}^\infty (1-zp^k)$, and
$$ \Gamma(z;p,q)=\prod_{j,k=0}^\infty
\frac{1-z^{-1}p^{j+1}q^{k+1}}{1-zp^jq^k}
$$
is the elliptic gamma function.
One can take $|t_a|>1$ and deform the contour of integration in such a way that it
continues to separate double geometric sequences of the poles of the integrand
$z_p^{in}=t_ap^jq^k,\, j,k\in\Z_{\geq 0},$ going to $z=0$ for $j,k\to\infty$,
from their reciprocals $z_p^{out}=t_a^{-1}p^{-j}q^{-k}$ going to $z=\infty$.
The right-hand side expression provides an analytical
continuation of the integral to all values of $t_a\in\CC^\times$ excluding the pole positions.

This is the most general known extension of the Gaussian integral as well as of the Euler beta integral.
At the very bottom it contains also the Newton's binomial theorem.
It is attached to the $E_6$ exceptional root system. For a special choice of parameters it can be reduced
to the Frenkel--Turaev sum \cite{die-spi:elliptic}.

The elliptic beta integral serves also as the simplest explicit example of the elliptic Fourier
transformation. With its help  the most general rank 1 solution of the Yang--Baxter equation has been
found and the corresponding solvable two-dimensional statistical mechanics system was constructed.
The most important application of this integral was found in quantum field theory.
In this context, formula \eqref{ellbeta} expresses the equality of superconformal indices of two
theories related by the Seiberg duality and describes the confinement phenomenon in the sector of 
BPS states of a simple four-dimensional $\mathcal{N}=1$ supersymmetric gauge field theory. 
A brief survey of these applications is given in \cite{spi:aar}.

Here we consider another point of view on the elliptic beta integral.
Namely, we try to represent it as the sum of residues of all integrand poles lying inside the unit circle.
Before passing to this problem, we describe the infinite very-well-poised elliptic hypergeometric series
\begin{eqnarray}
&&{}_{r+1}V_{r}(t_0;t_1,\ldots,t_{r-4};q,p):= \sum_{n=0}^\infty
\frac{\theta(t_0q^{2n};p)}{\theta(t_0;p)}\prod_{m=0}^{r-4}
\frac{\theta(t_m;p;q)_n}{\theta(w_m;p;q)_n}q^n,
\label{V-series}\end{eqnarray}
where $w_m=qt_0t_m^{-1}$. The parameters $t_m$ satisfy the balancing condition
\begin{equation}
\prod_{k=1}^{r-4}t_k=\nu t_0^{(r-5)/2}q^{(r-7)/2}, \quad \nu=\pm 1,
\quad \textrm{or}\quad \prod_{k=0}^{r-4}\frac{t_k}{w_k}=q^{-4}.
\label{balancing}\end{equation}
For odd $r$ the choice  $\nu=1$ is fixed by convention.
Here we use the theta function and the elliptic Pochhammer symbol
$$
\theta(a;p)=(a;p)_\infty(pa^{-1};p)_\infty, \qquad \theta(a;p;q)_k=\prod_{j=0}^{k-1}\theta(aq^j;p).
$$
When we talk on the radius of convergence of this infinite series we assume that the term $q^n$ in \eqref{V-series}
is replaced by $(qz)^n$ and the convergence is considered with respect to the values of the variable $z$.
The formal function $_{r+1}V_r$ is $p$-periodic in all its abelian parameters,
including $q$, i.e. it does not change under the transformations
$t_j\to p^{n_j}t_j,\, q\to qp^{m}$ respecting the
balancing condition, i.e. when $\sum_{k=1}^{r-4}n_j= \tfrac12 (r-5)n_0+\tfrac12 (r-7)m$.

Since the point $z=0$ is a non-isolated essential singularity of the integrand in \eqref{ellbeta},
computing the sum of residues is a delicate problem. We take as a model example a computation of the 
famous Askey-Wilson integral described in \cite{aw}. Its integrand depends on
four parameters, the basic variable $q$ and has infinitely many poles accumulating near the 
point $z=0$. Shrinking the contour of integration $\T$ to $C_\epsilon$ (a circle of an infinitesimally 
small radius $\epsilon$) one sees that the residues form four $_6\varphi_5$ series 
with a finite number of terms. For $|q|<1$ and some additional restriction on the parameters 
all four series are convergent and summable to the form of infinite products. Their emerging
combination can be represented as a unique infinite product yielding formally the value of the
integral. However, the convergence of the infinite sum of residues is not sufficient for 
such a claim, the fully 
rigorous proof needs an estimation of the error term coming from the integral over $C_\epsilon$
when $\epsilon\to 0$. Suppose an integral of a function with a non-isolated essential singularity 
at $z=0$ inside the contour of integration has the value $I$ and it is given by the convergent sum of
residues. Now we add to the integrand a function of $z$ having only
one isolated singularity at $z=0$ such that its Laurent expansion at this point
has a nonzero term $a/z$, $a\neq 0$ (actually, it sufficient to add just one pole $a/z$). 
Then the nature of singularity at $z=0$ does not change (it remains non-isolated), 
the infinite sums of residues of the integrand poles will be convergent,
but the limiting value of the integral will be $I+a$, acquiring in addition a finite value $a$. 
The error term in the second case will be convergent to $a$. Therefore the infinite sum of residues
represents the integral value only if the corresponding error term tends to zero.

Our situation is somewhat different from the one considered in \cite{aw}. 
We know already the value of the elliptic beta integral and we
would like to understand whether there is any kind of infinite series summation formula behind it. Now the
integrand in \eqref{ellbeta} has poles inside $\T$ lying at the points $z=t_ap^jq^k,\, a=1,\ldots,6,\, j,k\in\Z_{\geq 0}$.
Shrink $\T$ down to zero assuming that all the poles are simple, i.e. $t_b\neq t_ap^jq^k$, $b\neq a$.
The following formulas are useful for computation of residues:
$$
\theta(tp^j;p)=(-t)^{-j}p^{-\frac{j(j-1)}{2}}\theta(t;p), \quad j\in\Z,
$$
$$
\Gamma(tp^jq^k;p,q)=\frac{\theta(t;p;q)_k\theta(t;q;p)_j}{(-t)^{jk} q^{j\frac{k(k-1)}{2}}p^{k\frac{j(j-1)}{2}}}
\Gamma(t;p,q), \quad j,k\geq 0,
$$
$$
\Gamma(tp^{-j}q^{-k};p,q)=\frac{(-t)^{-jk-j-k} q^{(j+1)\frac{k(k+1)}{2}}p^{(k+1)\frac{j(j+1)}{2}}}
{\theta(qt^{-1};p;q)_k\theta(pt^{-1};q;p)_j} \Gamma(t;p,q), \quad j,k\geq 0.
$$

Then the infinite sum of residues yields (the details of computation are given in the Appendix)
\begin{eqnarray}\nonumber &&\makebox[-2em]{}
I(\underline{t};p,q)= \frac{1}{2} \sum_{a=1}^6\frac{\prod_{\ell=1,\neq a}^6\Gamma(t_\ell t_a^{\pm1};p,q)}{\Gamma(t_a^{-2};p,q)}
{}_{10}V_9(t_a^2; t_at_1,\ldots,\check t_a^2,\ldots, t_at_6;q,p)\Big|_{t_m\to \frac{t_m}{p}}
\\ && \makebox[4em]{} \times
{}_{10}V_9(t_a^2; t_at_1,\ldots,\check t_a^2,\ldots, t_at_6;p,q)\Big|_{t_m\to \frac{t_m}{q}},\quad m\neq a,
\label{sumres}\end{eqnarray}
where $\check t_a^2$ means the absence of this parameter and the replacement $t_m\to t_m/p$ or $t_m/q$
takes place only for one arbitrarily chosen value of $m\neq a$. 
Infinite $_{10}V_9$-series has the form
\begin{eqnarray} \nonumber &&\makebox[2em]{}
{}_{10}V_9(a; b,c,d,e,f;q,p)=
\\ && \makebox[-2em]{}
=\sum_{k=0}^\infty \frac{\theta(aq^{2k};p)}{\theta(a;p)}
\frac{\theta(a,b,c,d,e,f;p;q)_k\,q^k}{\theta(q, qa/b,qa/c,qa/d,qa/e,qa/f;p;q)_k}, \quad bcdef=qa^2.
\label{10V9}\end{eqnarray}
The parameters of ${}_{10}V_9$ series standing in \eqref{sumres} satisfy the needed balancing 
condition because 
one of the parameters (with the index $m$) is scaled either by $p^{-1}$ or $q^{-1}$. This keeps
the expression explicitly $p\leftrightarrow q$ symmetric. Such formal series representations
of the elliptic hypergeometric integrals was known to the author
shortly after the discovery of the latter integrals,
however their rigorous mathematical meaning for general values of parameters is absent until
the present time.

The key question is now whether the function $I(\underline{t};p,q)$,
as defined in \eqref{sumres}, is equal to
the right-hand side expression in \eqref{ellbeta}, representing thus a rigorously defined infinite
elliptic hypergeometric series identity? To answer this question, fisrt, one should analyze 
convergence of all ${}_{10}V_9$ series entering  \eqref{sumres} and, second, estimate the 
corresponding error term. We discuss here only the first point. The analysis of convergence 
of $_{r+1}V_r$ series
in \cite{KS} is not sufficient for answering this question. Below we describe a slightly
more general consideration, which still does not give the desired result.

Denote
$$
p=e^{2\pi i\tau},\, \quad \textup{Im}\,\tau>0, \quad q=e^{2\pi i\sigma},\, \quad \textup{Im}\,\sigma>0,
$$
and restrict the values of $q$ to
\begin{equation}
q=e^{2\pi i \chi(N+M\tau)}, \quad N, M\in\Z,\; (N,M)\neq (0,0), \quad M>0,
\label{q}\end{equation}
where $\chi$ is a real irrational number from the interval $[0,1]$. Equivalently, we can write
\begin{equation}
\sigma=\chi(N+M\tau) -Q, \quad Q\in\Z,\quad \mathrm{Im}\,\sigma=\chi M\mathrm{Im}\,\tau.
\label{pqrel1}\end{equation}
For convenience in the following, we fix $Q>0$.
As a result, for $n\in\Z_{\geq 0}$ the sequence $q^n\in e^{2\pi i\mathcal L}$, where $\mathcal L$ is
the line connecting the points $(0,0)$ and $(N,M\tau)$ in the lattice of points $\Z+\Z\tau$.

Next, it is necessary to demand that the parameters $t_m$ and $w_m=qt_0/t_m$
do not lie on the line $q^\RR$. Still, there remains a singularity in the denominator of
the series coefficients coming from the value $w_0=q$ which is approached arbitrarily close for $n\to\infty$.
In the work of Hardy and Littlewood \cite{HL}, similar singularity was investigated in the analysis of
convergence of a simple $q$-hypergeometric series when $|q|=1$, but $q$ is not a root of unity.
It was shown that the corresponding singularity is harmless, if
one constrains the values of $\chi$ to the irrational numbers, whose denominators $q_k$ of the
continued fraction approximants $\chi =\stackreb{\lim}{k\to \infty}p_k/q_k$ satisfy the constraint
\begin{equation}
\limsup_{k\to\infty} \frac{\log q_{k+1}}{q_k} = 0,
\label{HLcond}\end{equation}
which is true for almost all (i.e. of measure one) irrational numbers.

As shown in \cite{KS}, this picture holds for the general elliptic hypergeometric series as well
and the series \eqref{V-series} converge under the taken constraints on parameters for
sufficiently small $|z|$.
In order to determine the radius of convergence of the elliptic hypergeometric series,
the following integral has been computed exactly for $M, N, K\in\Z_{>0}$ and $(M,N)=1$:
\begin{eqnarray}\label{FNM} && \makebox[4em]{}
F_{KN,KM}(t):=\int_{0}^{1} \log |\theta(te^{2\pi i x K(N+M\tau)}; p)|d x
= \frac{\log^2|t|}{4\pi\mathrm{Im}\,\tau}
\\ &&
-\frac{(KM-1)\log|t|}{2}
+\frac{(KM-1)(2KM-1)\pi\mathrm{Im}\,\tau}{6}
-\frac{\mathrm{Im}\,\tau}{\pi|N+M\tau|^2}
\mathrm{Re}(\mathrm{Li}_2(\mu^{M})),
\nonumber\end{eqnarray}
where
\begin{equation}
\mu := \frac{t}{|t|}\exp\left(i (N+M\mathrm{Re}\, \tau)\frac{\log|t|}{M\mathrm{Im}\,\tau}\right)
=\exp\left(i\frac{\mathrm{Re} ((N+M\tau)\log \bar{t})}{M\mathrm{Im}\,\tau}\right),
\label{mu}\end{equation}
i.e. $|\mu|=1,$ and $\mathrm{Li}_2(x)$ is the standard Euler dilogarithm
$$
\mathrm{Li}_2(x)=\sum_{n=1}^\infty \frac{x^n}{n^2}, \quad x\in\CC, \; |x|\leqslant 1.
$$
There is a beautiful relation between the function $\mathrm{Li}_2(x)$ of the argument $|x|=1$
and the second Bernoulli polynomial \cite{aar},
\begin{equation}
\mathrm{Re}(\mathrm{Li}_2(x))=\pi^2B_2\left(\frac{\arg x}{2\pi}\right )= \frac{(\mathrm{arg}\, x)^2}{4}-\frac{\pi \mathrm{arg} \, x}{2}+\frac{\pi^2}{6}
\label{bern}\end{equation}
for $\mathrm{arg} \, x \in [0, 2\pi]$. This identity simplifies the expression \eqref{FNM} and converts it
into some quadratic polynomial.

We introduce the power counting term $z^n$ into series \eqref{V-series} and rewrite it in the form
$$
\sum_{n=0}^\infty c_n,\qquad c_n=\prod_{k=0}^{n-1} H(q^k)\, z^n,
$$
where $H(u)$ is an elliptic function of $u$ of the form
\begin{equation}
H(u)=q\frac{\theta(t_0q^2u^2;p)}{\theta(t_0u^2;p)}
\prod_{j=0}^{r-4} \frac{\theta(t_ju;p)}{\theta(qt_0u/t_j;p)}
\label{Hvwp}\end{equation}
with the parameters satisfying the well-poisedness $w_k=qt_0/t_k$ and balancing
$\prod_{k=0}^{r-4}t_k/w_k=q^{-4}$ conditions.
The function $H(u)$ is a $p$-elliptic function of $u$, $H(pu)=H(u)$, as well as of all its
independent abelian parameters, including $q$.

The modulus of the series coefficients can be represented in the form
$$
|c_n|
=\prod_{k=0}^{n-1}|H(q^k)|\, |z|^n
=\Big|z\exp \Big(\frac{1}{n} \sum_{k=0}^{n-1} \log |H(q^k)|\Big)\Big|^n.
$$
Under the taken constraints on the values of $\chi$, one can apply the Weyl
equidistribution theorem and write for $n\to\infty$ the leading asymptotics of these
coefficients
\begin{eqnarray} \nonumber &&
\stackreb{\lim}{n\to\infty}|c_n|=\stackreb{\lim}{n\to\infty}\exp\big(\sum_{k=0}^{n-1}\log|zH(q^k)|\big)
\propto \left|\frac{z}{r_c}\right|^n,\quad
\\ && \makebox[2em]{}
\log r_c^{-1}=\int_0^1\log|H(e^{2\pi i x(N+M\tau)})|dx,
\label{int_rc}\end{eqnarray}
i.e. $r_c$ is the radius of convergence of the corresponding series.
Evidently, the integral \eqref{int_rc} can be written as a combination
of the above integrals \eqref{FNM}. Namely, one has
$$
\log r_c^{-1}=\log|q|+ F_{2N,2M}(q^2t_0)-F_{2N,2M}(t_0)+\sum_{k=0}^{r-4}(F_{N,M}(t_k)-F_{N,M}(w_k)).
$$
Substituting formula \eqref{FNM} we obtain the following expression for the radius of convergence of $_{r+1}V_r$-series
\begin{eqnarray}\label{rc1} &&  \makebox[-1em]{}
\log r_{c}^{-1}=\frac{\pi\mathrm{Im}\,\tau}{|N+M\tau|^2}
\sum_{k=0}^{r-4} (\alpha_k(\alpha_k-1)-\beta_k(\beta_k-1)),
\\ &&  \makebox[1em]{}
\alpha_k:=\Biggl\{ \frac{\mathrm{Re}((N+M\tau)\log \bar{w}_k)}{2\pi\mathrm{Im}\,\tau}\Biggr\}, \quad
\beta_k:=\Biggl\{ \frac{\mathrm{Re}((N+M\tau)\log \bar{t}_k)}{2\pi\mathrm{Im}\,\tau}\Biggr\},
\label{rc1par}\end{eqnarray}
where $\{x\}$ means the fractional part of $x$. The expression for $r_c$ given in \cite{KS}
contains a wrong additional term $M\log|q|$ in the right-hand side of \eqref{rc1} (this error
occurred because of the use of combination $F_{N,M}(q^2t_0)-F_{N,M}(t_0)$ instead of the one given above).
Let us remind also that  $\alpha_0=0$ because $\mathrm{Re}((N+M\tau)\log \bar{q})=0$.

Note that any parameter $t_k$ or $w_k$ can be multiplied
by $q^h,\, h\in\mathbb{R},$ and  arbitrary integer power of $p$ and this does not influence
the radius of convergence. Indeed,
$$
\mathrm{Re}\big((N+M\tau)\log \bar q^h\bar t\big)=\mathrm{Re}\big((N+M\tau)\log \bar t\big)
$$
and
$$
\frac{\mathrm{Re}((N+M\tau)\log \bar{p}^m)}{2\pi\mathrm{Im}\,\tau}= -mN,
\quad \textup{and}\quad \{x-mN\}=\{x\}.
$$
Therefore whenever the series converges for a taken set of values of $t_k$, it will converge for
all  $t_kq^{h_k}p^{m_k}$, $m_k\in\Z$. This gives the non-intersecting spirals of convergence on the complex
plane of parameters $t_k$.

It is convenient to use the parametrization
\begin{equation}
t_j=q^{h_j} e^{\varphi_j  \frac{2\pi \mathrm{Im}\,\tau}{N+M\bar\tau}}, \qquad
w_j=\frac{qt_0}{t_j}=q^{\tilde{h}_j} e^{\tilde{\varphi}_j \frac{2\pi \mathrm{Im}\,\tau}{N+M\bar\tau}},
\quad h_j, \tilde{h}_j, \varphi_j, \tilde{\varphi}_j \in \mathbb{R},
\label{hphi}\end{equation}
where the variables $\varphi_j$ and $\tilde\varphi_j$ determine the distance of the poles
and zeros of the series from the line $q^\mathbb{R}$. Then,
\begin{equation}
\log r_{c}^{-1}=\frac{\pi\mathrm{Im}\,\tau}{|N+M\tau|^2} \sum_{k=0}^{r-4}  (\{\tilde{\varphi}_k\}-\{\varphi_k\})(\{\tilde{\varphi}_k\}+\{\varphi_k\}-1),
\label{rcvwp}\end{equation}
where $\tilde \varphi_0=0$ (as follows from the equality $w_0=q$) and
$$
\varphi_0=\varphi_k+\tilde \varphi_k, \; k=1,\ldots, r-4,\quad
\sum_{k=0}^{r-4} \varphi_k=\sum_{k=0}^{r-4}\tilde \varphi_k
$$
(the latter equality follows from the balancing condition).

If $\{\varphi_k\}=\varphi_k, \{\tilde\varphi_k\}=\tilde\varphi_k, \, k=0,\ldots, r-4$, then
all the terms in \eqref{rcvwp} vanish and $r_c=1$, i.e. the convergence of the original infinite
$_{r+1}V_r$ (i.e. without the damping factor $z^n$) very-well-poised elliptic hypergeometric
series \eqref{V-series}
remains under the question in this situation. However, in \cite{KS} a special choice of parameters
for the well poised elliptic hypergeometric series was suggested for which $r_c>1$
for arbitrary $r>2$:
\begin{eqnarray} \nonumber &&
\varphi_1=1+\frac{\varepsilon r}{2},\qquad \varphi_2=\ldots=\varphi_r =1-\varepsilon,
\quad \varepsilon >0,
\\ &&
\tilde{\varphi}_1=1-\varepsilon \Big(1-\frac{2}{r}+\frac{r}{2}\Big),\quad \tilde{\varphi}_2=\ldots=\tilde{\varphi}_r=1+\frac{2\varepsilon}{r},
\label{wppar}\end{eqnarray}
where
$$
\varepsilon = \frac{k+1}{\frac{r}{2} + \lambda},  \quad k=\left[\frac{r-2}{2}\right],
\quad 0< 2\lambda < 1-\frac{2}{r}.
$$
The same choice of parameters can be applied to the very-well-poised series
 $_{r+9}V_{r+8}$ when four parameters $t_{r+1}, \ldots, t_{r+4}$ are chosen in such a way
that the very-well-poised part of the coefficients $c_n$, $q^n\theta(t_0q^{2n};p)/\theta(t_0;p)$,
will be cancelled out.
Then, after taking $t_0=q$, one gets precisely the well poised series to which the convergence
constraints \eqref{wppar} apply.

Consider now the conditions of simultaneous convergence of all six $_{10}V_9(\ldots;q,p)$-series figuring
in \eqref{sumres}. Note that the parameters $t_k$ enter these series in a slightly different way than in \eqref{V-series}. We apply again the parametrization 
\begin{equation}
t_j=q^{h_j} e^{\varphi_j  \frac{2\pi \mathrm{Im}\,\tau}{N+M\bar\tau}}, \quad j=1,\ldots,6,
\quad \prod_{j=1}^6t_j=pq.
\label{hphires1}\end{equation}
One of the parameters $t_m$ is divided by $p$, which is irrelevant as mentioned above.
From the balancing condition we obtain the constraint $\sum_{j=1}^6\varphi_j=0$.
Then the radii of convergence of six different series take the form
\begin{equation}
-\log r_{c,q}^{(a)}=\frac{\pi\mathrm{Im}\,\tau}{|N+M\tau|^2} \sum_{k=1}^6
(\{\varphi_a-\varphi_k\}-\{\varphi_a+\varphi_k\})(\{\varphi_a-\varphi_k\}+\{\varphi_a+\varphi_k\}-1),
\label{rcall}\end{equation}
where $a=1,\ldots,6$. The subscript $``q"$ reminds that we deal with the elliptic
hypergeometric series with the basic variable $q$.
At the moment it is not clear whether there exists such a set
of variables $\varphi_k$ for which the sums on the right-hand side of equalities \eqref{rcall} all
are negative, which would mean that all $r_{c,q}^{(a)}>1$.

The derived result is not sufficient for analyzing simultaneous convergence of
$_{r+1}V_r$-series with different modular  parameters $p$ and $q$ in \eqref{sumres}.
To be able to treat them uniformly we have to consider convergence of the
original $_{r+1}V_r$-series for a slightly more general choice of the parameter $q$.
Namely, let us set
\begin{equation}
q=\epsilon e^{2\pi i \chi (N+M\tau)},\quad  \epsilon= e^{2\pi i L/K}, \quad \epsilon^K=1,
\label{q_gen}\end{equation}
where $N, M, L, K \in \Z,\, (N,M)=(L,K)=1.$

Using the representation $n=mK+j$, $j=0,1,\ldots, K-1$, $m\in\Z_{\geq 0}$, we come to
the Weyl equidistribution theorem in the following form
\begin{eqnarray}\nonumber &&
\stackreb{\lim}{n\to\infty}\frac{1}{n}\sum_{k=0}^{n-1}\log|H(q^k)|
=\stackreb{\lim}{m\to\infty}\sum_{j=0}^{K-1}\frac{1}{mK+j}\sum_{\ell=0}^{m}
\log|H(\epsilon^j e^{2\pi i (\ell K+j)\chi(N+M\tau)})|
\\ && \makebox[3em]{}
=\frac{1}{K}\sum_{j=0}^{K-1}\int_0^1\log|H(q^j e^{2\pi i xK(N+M\tau)})|dx.
\label{Weyl_gen}\end{eqnarray}
Now we use formula \eqref{FNM} and see that the expression \eqref{rc1} is modified
as follows
\begin{eqnarray}\label{rc2} &&  \makebox[-1em]{}
\log r_{c}^{-1}=\frac{\pi\mathrm{Im}\,\tau}{|N+M\tau|^2}
\sum_{j=0}^{K-1}\sum_{k=0}^{r-4} (\alpha_k^{(j)}(\alpha_k^{(j)}-1)-\beta_k^{(j)}(\beta_k^{(j)}-1)),
\\ &&  \makebox[1em]{}
\alpha_k^{(j)}:=\Biggl\{ \frac{\mathrm{Re}((N+M\tau)
\log \overline{\epsilon^jw_k})}{2\pi\mathrm{Im}\,\tau}\Biggr\}, \quad
\beta_k^{(j)}:=\Biggl\{ \frac{\mathrm{Re}((N+M\tau)
\log \overline{\epsilon^jt_k})}{2\pi\mathrm{Im}\,\tau}\Biggr\}.
\label{rc2par}\end{eqnarray}
Using the same representation of parameters \eqref{hphi}, we obtain
\begin{eqnarray}\nonumber && \makebox[4em]{}
\log r_{c,q}^{-1}=\frac{\pi\mathrm{Im}\,\tau}{|N+M\tau|^2}
\\ && \times
 \sum_{j=0}^{K-1}\sum_{k=0}^{r-4}  (\{\tfrac{LM}{K}j+\tilde{\varphi}_k\}-\{\tfrac{LM}{K}j+\varphi_k\})(\{\tfrac{LM}{K}j+\tilde{\varphi}_k\}
+\{\tfrac{LM}{K}j+\varphi_k\}-1),
\label{rcvwp2}\end{eqnarray}
with the same restrictions on the variables $\varphi_k$ and $\tilde \varphi_k$ as before.

Suppose that the series $_{r+1}V_r$ is convergent. Then it defines
a doubly periodic function of parameters $g_k$, $t_k=e^{2\pi i g_k}$,
with natural boundaries in the fundamental
parallelogram of periods. These boundaries are formed by the sequences
$$
Lj/K + n\chi(N+M\tau) \mod 1, \tau,
$$
which fill for $n\to\infty$ not one line
on complex plane of abelian parameters stretching from the zero point to $N+M\tau$,
but $K$ parallel lines starting from rational points $0, 1/K, \ldots, (K-1)/K$.

As to the six $_{10}V_9(\ldots;q,p)$ series entering \eqref{sumres}, their radii of convergence are
given by the following natural generalization of formula \eqref{rcall},
\begin{eqnarray}\nonumber &&
-\log r_{c,q}^{(a)}=\frac{\pi\mathrm{Im}\,\tau}{|N+M\tau|^2}
 \sum_{j=0}^{K-1}\sum_{k=1}^6 (\{\tfrac{LM}{K}j+\varphi_a-\varphi_k\}
 -\{\tfrac{LM}{K}j+\varphi_a+\varphi_k\})
 \\ && \makebox[4em]{}
 \times (\{\tfrac{LM}{K}j+\varphi_a-\varphi_k\}
 +\{\tfrac{LM}{K}j+\varphi_a+\varphi_k\}-1),
\label{rcallgen}\end{eqnarray}
where $a=1,\ldots,6$ and $\sum_{k=1}^6\varphi_k=0$.

Now we have to investigate convergence of the $_{10}V_9(\ldots;p,q)$-series in \eqref{sumres},
or of the series \eqref{V-series} with permuted $p$ and $q$.
For that we should invert relation \eqref{q} and express the base variable $p=e^{2\pi i\tau}$
in terms of $q$. From \eqref{pqrel1} it follows that
$$
\tau= \frac{\sigma+Q}{M\chi} -\frac{N}{M}, \quad Q, M, N\in \Z_{>0}.
$$
Let us write
$$
\frac{1}{M\chi}=S+\eta,\quad S=\left[\tfrac{1}{M\chi}\right]\in\Z_{\geq 0}, \quad 0<\eta<1,
$$
where $[x]$ is the integer part of $x$.
As mentioned above, due to the total ellipticity, the $_{r+1}V_r$-series under consideration does not change when $p$ is replaced  by $pq^S$. Therefore we can discard $S$ and write
\begin{equation}
p=e^{2\pi i\eta(\sigma+Q)}.
\label{pqrel2}\end{equation}
This situation is equivalent to the previously considered one with $\tau, N, M, L, K$ replaced
by $\sigma, Q, 1, -N, M$, respectively. In particular, $\epsilon$ is replaced by
$\omega=e^{-2\pi i\frac{N}{M}}$. It is necessary also to impose a constraint on the
values of $\eta$ --- the denominators of its continued fraction approximants should satisfy the same
constraint \eqref{HLcond}, which is true for the measure one set of numbers. The latter means
that the two constraints, on $\chi$ and $\eta$, are simultaneously satisfied for a measure
one set of numbers. As a result, we can write immediately the expression
for the radius of convergence
\begin{eqnarray}\label{rc3} &&  \makebox[-1em]{}
\log r_{c}^{-1}=\frac{\pi\mathrm{Im}\,\sigma}{|Q+\sigma|^2}
\sum_{j=0}^{M-1}\sum_{k=0}^{r-4} (\gamma_k^{(j)}(\gamma_k^{(j)}-1)-\delta_k^{(j)}(\delta_k^{(j)}-1)),
\\ &&  \makebox[1em]{}
\gamma_k^{(j)}:=\Biggl\{ \frac{\mathrm{Re}((Q+\sigma)
\log \overline{\omega^jw_k})}{2\pi\mathrm{Im}\,\sigma}\Biggr\}, \quad
\delta_k^{(j)}:=\Biggl\{ \frac{\mathrm{Re}((Q+\sigma)
\log \overline{\omega^jt_k})}{2\pi\mathrm{Im}\,\sigma}\Biggr\}.
\label{rc3par}\end{eqnarray}
Apply now a new parametrization
\begin{equation}
t_j=p^{f_j} e^{\xi_j  \frac{2\pi \mathrm{Im}\,\sigma}{Q+\bar\sigma}}, \quad
w_j=p^{\tilde f_j} e^{\tilde\xi_j  \frac{2\pi \mathrm{Im}\,\sigma}{Q+\bar\sigma}},
\quad f_j, \tilde f_j, \xi_j, \tilde \xi_j \in \mathbb{R}.
\label{hphires2}\end{equation}
Then the radius of convergence of the  $_{r+1}V_r(\ldots;p,q)$-series takes the form
\begin{eqnarray}\nonumber && \makebox[4em]{}
\log r_{c,p}^{-1}=\frac{\pi\mathrm{Im}\,\tau}{|Q+\sigma|^2}
\\ && \times
 \sum_{j=0}^{M-1}\sum_{k=0}^{r-4}  (\{\tilde{\xi}_k-\tfrac{N}{M}j\}-\{\xi_k-\tfrac{N}{M}j\})
 (\{\tilde{\xi}_k-\tfrac{N}{M}j\} +\{\xi_k-\tfrac{N}{M}j\}-1),
\label{rcvwp3}\end{eqnarray}
where $\tilde\xi_0=0$, $\xi_0=\tilde\xi_j+\xi_j,\, j=1,\ldots, r-4$, and
$\sum_{j=0}^{r-4} (\tilde\xi_j-\xi_j)=0$.

Analogously, the radii of convergence of six $_{10}V_9(\ldots;p,q)$-series entering
the equality \eqref{sumres} have the form
\begin{eqnarray}\label{rcallp} && \makebox[2em]{}
-\log r_{c,p}^{(a)}=\frac{\pi\mathrm{Im}\,\sigma}{|Q+\sigma|^2}
\\ && \times
\sum_{j=0}^{M-1}\sum_{k=1}^6
(\{\xi_a-\xi_k-\tfrac{N}{M}j\}-\{\xi_a+\xi_k-\tfrac{N}{M}j\})(\{\xi_a-\xi_k-\tfrac{N}{M}j\}
+\{\xi_a+\xi_k-\tfrac{N}{M}j\}-1),
\nonumber\end{eqnarray}
where $a=1,\ldots,6$ and $\sum_{k=1}^6\xi_k=0$. We conclude that, in order to partially 
justify the representation of the elliptic beta integral as a quadratic combination of 12 infinite
$_{10}V_9$ very-well-poised elliptic hypergeometric series, it is necessary
to find the values of variables $\varphi_k$ and $\xi_k$, emerging from the parametrizations
$$
t_j=q^{h_j} e^{\varphi_j  \frac{2\pi \mathrm{Im}\,\tau}{N+M\bar\tau}}
=p^{f_j} e^{\xi_j  \frac{2\pi \mathrm{Im}\,\sigma}{Q+\bar\sigma}},\quad j=1,\ldots, 6,
$$
such that the right-hand side expressions of all 12 equalities \eqref{rcallgen} and \eqref{rcallp}
are negative.

\section{Asymptotics of the Frenkel--Turaev sum}

The Frenkel--Turaev sum, derived in \cite{FT}, provides a closed form expression for the
terminating $_{10}V_9$ very-well-poised elliptic hypergeometric series
\begin{eqnarray}\nonumber &&
{}_{10}V_9(t_0;t_1,\ldots,t_5;q,p)=\sum_{k=0}^n\frac{\theta(t_0q^{2k};p)}{\theta(t_0;p)}
\prod_{m=0}^5 \frac{\theta(t_m;p;q)_k} {\theta(qt_0/t_m;p;q)_k}\, q^k
\\ && \makebox[1em]{}
=\frac{\theta(qt_0,\frac{qt_0}{t_1t_2},
\frac{qt_0}{t_1t_3},\frac{qt_0}{t_2t_3};p;q)_n }
{\theta(\frac{qt_0}{t_1t_2t_3},\frac{qt_0}{t_3},
\frac{qt_0}{t_2},\frac{qt_0}{t_1};p;q)_n}
=: S_n(p), \quad t_4=q^{-n}, \quad \prod_{m=1}^5t_m=qt_0^2.
\label{FTsum}\end{eqnarray}
It can be checked that both sides of this identity are $p$-elliptic functions of all parameters $t_j$ and $q$,
i.e. we have a particular elliptic functions identity.

For fixed $z$ and $p\to 0$ one has $\theta(z;0)=1-z$ and the elliptic Pochhammer symbol degenerates to the $q$-shifted factorial
$(z;q)_n=\prod_{k=0}^{n-1}(1-zq^k)$. Therefore, for fixed parameters $t_k$ this limit degenerates ${}_{10}V_9$-series to the very-well-poised balanced
${}_{8}W_7$ $q$-hypergeometric series and one gets the Jackson sum \cite{GR}:
\begin{eqnarray} \nonumber &&
{}_{8}W_7(t_0; t_1,...,t_5; q, q)=\sum_{k=0}^n\frac{1-t_0q^{2k}}{1-t_0}
\prod_{m=0}^5 \frac{(t_m;q)_k} {(qt_0/t_m;q)_k}\, q^k
\\ && \makebox[2em]{}
=\frac{(qt_0,\frac{qt_0}{t_1t_2},
\frac{qt_0}{t_1t_3},\frac{qt_0}{t_2t_3};q)_n }
{(\frac{qt_0}{t_1t_2t_3},\frac{qt_0}{t_3},
\frac{qt_0}{t_2},\frac{qt_0}{t_1};q)_n}.
\label{Jack_sum}\end{eqnarray}
Assume now that $|q|<1$ and take the limit $n\to\infty$. This leads to the summation formula
for an infinite ${}_6W_5$-series appearing in the sum of residues representation 
of the Askey-Wilson integral
\begin{equation}
{}_{8}W_7
\stackreb{=}{n\to\infty} {}_{6}W_5(t_0;t_1,t_2,t_3;\frac{qt_0}{t_1t_2t_3})=
\frac{(qt_0,\frac{qt_0}{t_1t_2},\frac{qt_0}{t_1t_3},\frac{qt_0}{t_2t_3};q)_\infty }
{(\frac{qt_0}{t_1t_2t_3},\frac{qt_0}{t_3},\frac{qt_0}{t_2},\frac{qt_0}{t_1};q)_\infty}.
\label{8W7}\end{equation}

As to the Frenkel--Turaev sum, such a limit has a qualitatively different character.
 Note that on the right-hand side of \eqref{FTsum} we have the well poised
pattern of the ratios of theta functions---the products of parameters in the arguments of
theta functions in the numerator and denominator
are equal to $q^2t_0^2/t_1t_2t_3$. Therefore we may apply the analysis of the convergence of
well poised elliptic hypergeometric series performed in \cite{KS} to the present situation.
For fixed $t_k$ and $p\neq0$ consider $\lim_{n\to\infty}|S_n|$ for $q$ of the form \eqref{q}.
This means that we are interested in the asymptotics of the modulus of a particular elliptic function
when its order goes to
infinity in a very special way. For arbitrary irrational $\chi\in[0,1]$, we obtain
$$
|S_n(p)|=e^{\sum_{k=1}^n\log|H(q^k)|}\stackreb{=}{n\to\infty}\exp\big(n\int_0^1\log|H(e^{2\pi ix(N+M\tau)})|dx+o(n)\big),
$$
where
$$
H(u)=\frac{\theta(ut_0,\frac{ut_0}{t_1t_2},
\frac{ut_0}{t_1t_3},\frac{ut_0}{t_2t_3};p) }
{\theta(\frac{ut_0}{t_1t_2t_3},\frac{ut_0}{t_3},
\frac{ut_0}{t_2},\frac{ut_0}{t_1};p)}.
$$

Applying the results of \cite{KS} and using, as in the previous section, the parametrization
$t_j=q^{h_j} e^{\varphi_j  \frac{2\pi \mathrm{Im}\,\tau}{N+M\bar\tau}}$, we come to the integral evaluation
\begin{eqnarray}\nonumber && 
I_{FT}=\int_0^1\log|H(e^{2\pi ix(N+M\tau)})|dx=\kappa c,\qquad
\kappa= \frac{\pi\mathrm{Im} \,\tau}{|N+M\tau|^2}>0,
\\ \label{FTint}&&
c=(\{\varphi_0-\sum_{k=1}^3\varphi_k\}-\{\varphi_0\})(\{\varphi_0-\sum_{k=1}^3\varphi_k\}+\{\varphi_0\}-1)
\\ && \makebox[2em]{}
+\sum_{j=1}^3 (\{\varphi_0-\varphi_j\}-\{\varphi_0+\varphi_j-\sum_{k=1}^3\varphi_k\})
(\{\varphi_0-\varphi_j\}+\{\varphi_0+\varphi_j-\sum_{k=1}^3\varphi_k\}-1).
\nonumber \end{eqnarray}
There are three qualitatively different situations emerging for different values of
unconstrained $\varphi_k$-variables.
\begin{itemize}

\item
If all $\{x\}=x$, then $I_{FT}=0$ and $|S_n|\stackreb{=}{n\to\infty} e^{o(n)}$, i.e. we can only say that the exponential rate of grows or vanishing of the sum is absent. For a better estimate of the asymptotics one should know the behaviour of the error term in the Weyl equidistribution theorem.

\medskip

\item A convergence example. If $\varphi_0=3.5,\, \varphi_{1,2,3}=1.2$, then
\begin{eqnarray*} &&
c=(\{-0.1\}-\{3.5\})(\{-0.1\}+\{3.5\}-1)+3(\{2.3\}-\{1.1\})(\{2.3\}+\{1.1\}-1)
\\ && \makebox[2em]{}
=(0,4)^2- (0,6)^2=-0.2
\end{eqnarray*}
and $|S_n|\stackreb{=}{n\to\infty} e^{-0.2\kappa n}\to 0$.
The Frenkel--Turaev sum vanishes exponentially fast.

\medskip

\item  A divergence example. If $\varphi_0=0.1,\, \varphi_{1,2,3}=1.2$, then
\begin{eqnarray*} &&
c=(\{-3.5\}-\{0.1\})(\{-3.5\}+\{0.1\}-1)+3(\{-1.1\}-\{-2.3\})(\{-1.1\}+\{-2.3\}-1)
\\ && \makebox[2em]{}
=-(0,4)^2 + (0,6)^2 = 0.2
\end{eqnarray*}
and $|S_n|\stackreb{=}{n\to\infty} e^{0.2\kappa n}\to  \infty$.
The Frenkel--Turaev sum blows up exponentially fast.

\end{itemize}

We stress that these estimates are valid for arbitrary irrational $\chi$ despite of the presence of the singular $\theta(q;p;q)_n$ factor
in the denominator of the series terms. This happens because the termination condition
violates the constraint for parameters not to lie on the line $q^\RR$ which was used
in \cite{KS} in the convergence analysis of elliptic hypergeometric series.

Let us choose $M=0$, which reduces the value of $N$ to $N=1$. Then, taking the limit $p\to0$,
we can estimate the $n\to\infty$ limit of the Jackson sum \eqref{Jack_sum}
when $q$ lies on the unit circle, $q=e^{2\pi i\chi}$ with arbitrary irrational $\chi$.
The situation is now much simpler than before. According to \cite{HL,KS}, in the limit
$p\to 0$, one has
\begin{equation}\nonumber
F_{1,0}(t)=\int_{0}^{1} \log |1-te^{2\pi ix}|dx
= \begin{cases}
      0, & \text{if}\ |t|\leq 1 \\
      \log|t|, & \text{if}\ |t|>1
    \end{cases}
\end{equation}
As a result, we have
$$
|S_n(0)|=e^{cn+o(n)}, \quad c= \log \Big|\frac{(t_0)_{>1}
(t_0/t_1t_2)_{>1}(t_0/t_1t_3)_{>1}(t_0/t_2t_3)_{>1}}
{(t_0/t_1t_2t_3)_{>1}(t_0/t_3)_{>1}(t_0/t_2)_{>1}(t_0/t_1)_{>1}}\Big|,
$$
where $(t)_{>1}=1$ for $|t|\leq1$ and $(t)_{>1}=t$ for $|t|>1$.
Let us give examples of three possible regimes of behaviour.
If $|t_j|<1$, $j=1,2,3$, and $|t_0|<|t_1t_2t_3|$, then $c=0$.
If $1<|t_0|<|t_j|,$ $j=1,2,3$, then $c=\log|t_0|>0$ and $S_n(0)$
blows up exponentially fast. If $|t_0|<1<|t_0/t_1t_2t_3|$ and all other
ratios of parameters of interest are smaller than 1, then $c=-\log|t_0/t_1t_2t_3|<0$ and
$S_n(0)$ vanishes exponentially fast.
In the elliptic case, it is desirable to give somewhat similar constructive
description of the domains of $\varphi_k$ characterizing all three possible regimes,
not to limit to particular numerical examples.

\section{Conclusion}

Suppose that there exist such $\varphi_j$ and $\xi_j$ that all elliptic hypergeometric series
in the expression \eqref{sumres} converge and define analytical function $I(\underline{t};p,q)$
of parameters beyond the natural boundaries.
Then the question arises, how  their particular combination
of interest defines the meromorphic function of parameters described by the
product of elliptic gamma functions on the right-hand side of equality \eqref{ellbeta},
i.e. how the natural boundaries would disappear? A related question is whether there exists
a simpler elliptic hypergeometric series identity, a kind of elliptic extension
of the Jackson sum for non-terminating $_8\varphi_7$ $q$-hypergeometric series \cite{GR},
which would prove the residues sum representation of the elliptic beta integral
from scratch? Another open problem is the analysis of convergence of infinite elliptic
hypergeometric series for the general basic parameter $q$, which requires a two-dimensional
generalization of the Weyl equidistribution theorem and corresponding extension
of the Hardy-Littlewood convergence criterion \cite{HL,KS}.

\appendix

\section{Computation of residues}

On the left-hand side of \eqref{ellbeta} we shrink the integration contour, apply the
Cauchy theorem and pick up the residues. This yields
$$
2I(\underline{t};p,q)=\frac{(p;p)_\infty(q;q)_\infty}{2\pi i}\sum_{a=1}^6 \sum_{j,k=0}^\infty 2\pi i\stackreb{\lim}{z\to t_ap^jq^k}(1-\frac{t_ap^jq^k}{z})
\frac{\prod_{\ell=1}^6\Gamma(t_\ell z^{{\pm 1}};p,q)}{\Gamma(z^{\pm 2};p,q)}
$$

$$
=(p;p)_\infty(q;q)_\infty\sum_{a=1}^6 \sum_{j,k=0}^\infty \stackreb{\lim}{z\to t_ap^jq^k}(1-\frac{t_ap^jq^k}{z})\Gamma(\frac{t_a}{z}p^jq^k;p,q)
 \frac{\Gamma(\frac{t_a}{z};p,q)}{\Gamma(\frac{t_a}{z}p^jq^k;p,q)}
$$

$$ \times
\frac{\Gamma(t_a^2p^{j}q^{k};p,q)}{\Gamma(t_a^2p^{2j}q^{2k};p,q)\Gamma(t_a^{-2}p^{-2j}q^{-2k};p,q)}
\prod_{\ell=1,\neq a}^6\Gamma(t_\ell t_a p^{j}q^{k};p,q)\Gamma(\frac{t_\ell}{t_a} p^{-j}q^{-k};p,q)
$$

$$
= \sum_{j,k=0}^\infty \sum_{a=1}^6\stackreb{\lim}{\varepsilon \to 1}
\frac{\Gamma(\varepsilon p^{-j}q^{-k};p,q)}{\Gamma(\varepsilon;p,q)}
\frac{\prod_{\ell=1,\neq a}^6\Gamma(t_\ell t_a^{\pm1};p,q)}{\Gamma(t_a^{-2};p,q)}
$$

$$\times
\frac{\Gamma(t_a^2p^{j}q^{k};p,q)}{\Gamma(t_a^2p^{2j}q^{2k};p,q)}
\frac{\Gamma(t_a^{-2};p,q)}{\Gamma(t_a^{-2}p^{-2j}q^{-2k};p,q)}
\prod_{\ell=1,\neq a}^6\frac{\Gamma(t_\ell t_a p^{j}q^{k};p,q)}{\Gamma(t_\ell t_a;p,q)}\frac{\Gamma(\frac{t_\ell}{t_a} p^{-j}q^{-k};p,q)}{\Gamma(\frac{t_\ell}{t_a} ;p,q)}
$$

$$
= \sum_{a=1}^6\frac{\prod_{\ell=1,\neq a}^6\Gamma(t_\ell t_a^{\pm1};p,q)}{\Gamma(t_a^{-2};p,q)} \sum_{j,k=0}^\infty \frac{(-1)^{-jk-j-k} q^{(j+1)\frac{k(k+1)}{2}}p^{(k+1)\frac{j(j+1)}{2}}}
{\theta(q;p;q)_k\theta(p;q;p)_j}
$$

$$\times
\prod_{\ell=1,\neq a}^6\frac{\theta(t_\ell t_a;p;q)_k\theta(t_\ell t_a;q;p)_j}{(-t_\ell t_a)^{jk} q^{j\frac{k(k-1)}{2}}p^{k\frac{j(j-1)}{2}}} \frac{(-\frac{t_\ell}{t_a})^{-jk-j-k} q^{(j+1)\frac{k(k+1)}{2}}p^{(k+1)\frac{j(j+1)}{2}}}
{\theta(q\frac{t_a}{t_\ell};p;q)_k\theta(p\frac{t_a}{t_\ell};q;p)_j}
$$

$$\times
\frac{\theta(t_a^{2};p;q)_k\theta(t_a^{2};q;p)_j}{(-t_a^{2})^{jk} q^{j\frac{k(k-1)}{2}}p^{k\frac{j(j-1)}{2}}}
\frac{(-t_a^{2})^{4jk} q^{2j\frac{2k(2k-1)}{2}}p^{2k\frac{2j(2j-1)}{2}}}{\theta(t_a^{2};p;q)_{2k}\theta(t_a^{2};q;p)_{2j}}
$$

$$ \times
\frac{\theta(qt_a^{2};p;q)_{2k}\theta(pt_a^{2};q;p)_{2j}}
{(-t_a^{-2})^{-4jk-2j-2k} q^{(2j+1)\frac{2k(2k+1)}{2}}p^{(2k+1)\frac{2j(2j+1)}{2}}}
$$

$$
= \frac{1}{2} \sum_{a=1}^6\frac{\prod_{\ell=1,\neq a}^6\Gamma(t_\ell t_a^{\pm1};p,q)}{\Gamma(t_a^{-2};p,q)} \sum_{j,k=0}^\infty \frac{\theta(t_a^2q^{2k};p)}{\theta(t_a^2;p)} \frac{\theta(t_a^2p^{2j};q)}{\theta(t_a^2;q)}
\prod_{\ell=1}^6\frac{\theta(t_\ell t_a;p;q)_k\theta(t_\ell t_a;q;p)_j}
{\theta(q\frac{t_a}{t_\ell};p;q)_k\theta(p\frac{t_a}{t_\ell};q;p)_j}
$$

$$\times
\prod_{\ell=1}^6 \frac{(-\frac{t_\ell}{t_a})^{-jk-j-k} q^{(j+1)\frac{k(k+1)}{2}}p^{(k+1)\frac{j(j+1)}{2}}}
{(-t_\ell t_a)^{jk} q^{j\frac{k(k-1)}{2}}p^{k\frac{j(j-1)}{2}}}
\frac{(-t_a^{2})^{4jk} q^{2j\frac{2k(2k-1)}{2}}p^{2k\frac{2j(2j-1)}{2}}}
{(-t_a^{-2})^{-4jk-2j-2k} q^{(2j+1)\frac{2k(2k+1)}{2}}p^{(2k+1)\frac{2j(2j+1)}{2}}}
$$

$$
= \sum_{a=1}^6\frac{\prod_{\ell=1,\neq a}^6\Gamma(t_\ell t_a^{\pm1};p,q)}{\Gamma(t_a^{-2};p,q)} \sum_{j,k=0}^\infty \frac{\theta(t_a^2q^{2k};p)}{\theta(t_a^2;p)} \frac{\theta(t_a^2p^{2j};q)}{\theta(t_a^2;q)}
$$

$$ \times
\prod_{\ell=1}^6\frac{\theta(t_\ell t_a;p;q)_k\theta(t_\ell t_a;q;p)_j}
{\theta(q\frac{t_a}{t_\ell};p;q)_k\theta(p\frac{t_a}{t_\ell};q;p)_j}
q^{k(k+1)-j}p^{j(j+1)-k}t_a^{2j+2k}
$$

$$
= \sum_{a=1}^6\frac{\prod_{\ell=1,\neq a}^6\Gamma(t_\ell t_a^{\pm1};p,q)}{\Gamma(t_a^{-2};p,q)} \sum_{j,k=0}^\infty \frac{\theta(t_a^2q^{2k};p)}{\theta(t_a^2;p)} \frac{\theta(t_a^2p^{2j};q)}{\theta(t_a^2;q)}
$$

$$\times
\prod_{\ell=1,\neq m}^6\frac{\theta(t_\ell t_a;p;q)_k\theta(t_\ell t_a;q;p)_j}
{\theta(q\frac{t_a}{t_\ell};p;q)_k\theta(p\frac{t_a}{t_\ell};q;p)_j}
\frac{\theta(\frac{t_m}{p} t_a;p;q)_k\theta(\frac{t_m}{q} t_a;q;p)_j}
{\theta(pq\frac{t_a}{t_m};p;q)_k\theta(pq\frac{t_a}{t_m};q;p)_j}q^{k}p^{j}
$$

$$
= \sum_{a=1}^6\frac{\prod_{\ell=1,\neq a}^6\Gamma(t_\ell t_a^{\pm1};p,q)}{\Gamma(t_a^{-2};p,q)}
{}_{10}V_9(t_a^2; t_at_1,\ldots,\check t_a^2,\ldots, t_at_6;q,p)\Big|_{t_m\to \frac{t_m}{p}},
$$

$$ \times
{}_{10}V_9(t_a^2; t_at_1,\ldots,\check t_a^2,\ldots, t_at_6;p,q)
\Big|_{t_m\to \frac{t_m}{q}},\quad m\neq a.
$$
Here $\check t_a^2$ means that this parameter is absent in the arguments of $_{10}V_9$-functions \eqref{10V9}. The permutational symmetry in parameters is present, because the scalings
of one of the parameters $t_m/p$ and $t_m/q$ are equivalent for all $m\neq a$.

\medskip

{\bf Acknowledgments.}
The author is indebted to V.K. Beloshapka and S.M. Khoroshkin for useful discussions and 
to the referee for constructive remarks.
This study has been partially supported by the Russian Science Foundation (grant 24-21-00466).

\bibliographystyle{unsrt}

\end{document}